\documentclass[final]{article}

\usepackage{amsmath,amssymb,amsthm,amsfonts}
\usepackage{epsfig}

\usepackage[color,notref,notcite]{showkeys}
\definecolor{labelkey}{rgb}{0.6,0,1}

\newcommand{\R}{\mathbb{R}}
\renewcommand{\S}{\mathbb{S}_d}

\newcommand{\N}{\mathbb{N}}
\newcommand{\eps}{\varepsilon}
\newcommand{\un}{\mathbf{1}}

\def\xb{\bar{x}}
\def\yb{\bar{y}}
\def\tb{\bar{t}}
\def\ou{\overline u}
\def\uu{\underline u}

\def\oF{\overline F}
\def\uF{\underline F}
\def\ue{u^{\varepsilon}}
\def\ze{z^{\varepsilon}}
\def\Fe{F_{\varepsilon}}
\def\Fa{F_{\alpha}}
\def\limssup{\mathop{\rm limsup\!^*}}
\def\limiinf{\mathop{\rm liminf_*}}
\def\hrl{{\displaystyle{\mathop{\scriptstyle{y\to
x}}_{\varepsilon \to 0}}}}

\newtheorem{defi}{Definition}
\newtheorem{lem}{Lemma}
\newtheorem{prop}{Proposition}
\newtheorem{thm}{Theorem}
\newtheorem{cor}{Corollary}

\theoremstyle{remark}
\newtheorem{rem}{Remark}
\newtheorem*{rem*}{Remark}
\newtheorem{rems}{Remarks}
\theoremstyle{remark}
\newtheorem{ex}{Example}
\theoremstyle{remark}

\begin{document}

\title{Second-Order Elliptic Integro-Differential Equations:\\
Viscosity Solutions' Theory Revisited}

\author{Guy Barles\footnote{Laboratoire de Math\'ematiques et Physique Th\'eorique
CNRS UMR 6083, F\'ed\'eration Denis Poisson,
Universit\'e Fran\c{c}ois Rabelais, Parc de Grandmont,
37200 Tours, France, \texttt{barles@lmpt.univ-tours.fr}}\ \ and  
Cyril Imbert\footnote{Polytech'Montpellier \& 
Institut  de math\'ematiques et de mod\'elisation de Montpellier, UMR CNRS 5149,
Universit\'e Montpellier~II, CC 051, Place E. Bataillon, 34 095
Montpellier cedex 5, France, \texttt{imbert@math.univ-montp2.fr}}}

\maketitle

%%%%%%%%%%%%%%%%%%%%%%%%%%%%%%%%%%%%%%%%%%%%%%%%%%%%%%%%%%%%%%%%%%%%%%%%%%%%%
\begin{quote} \footnotesize
\noindent \textsc{Abstract.} The aim of this work is to revisit viscosity solutions' theory for second-order elliptic integro-differential equations and to provide a general framework which takes into account solutions with arbitrary growth at infinity. Our main contribution is a new Jensen-Ishii's Lemma for integro-differential equations, which is stated for solutions with no restriction on their growth at infinity. The proof of this result, which is of course a key ingredient to prove comparison principles, relies on a new definition of viscosity solution for integro-differential equation (equivalent to the two classical ones) which combines the approach with test-functions and sub-superjets.
\end{quote}
%%%%%%%%%%%%%%%%%%%%%%%%%%%%%%%%%%%%%%%%%%%%%%%%%%%%%%%%%%%%%%%%%%%%%%%%%%%%%%%
\vspace{5mm}

\medskip


\section*{Introduction}

In this paper, we revisit viscosity solutions' theory for second-order elliptic integro-differential equations. We first present in a general framework three (equivalent) definitions of viscosity solutions, each of them having its own utility; in particular, one of them is new. We also present stability results and we discuss comparison principles on a model example.  
\medskip

The core of the paper lies in an analogue of the celebrated Jensen-Ishii's Lemma in  the framework of second-order elliptic integro-differential equations. This nonsmooth analysis lemma is the keystone of the proofs of comparison principles in 
viscosity solution theory for \emph{local} second-order fully nonlinear elliptic equations but, because of some particular features of \emph{nonlocal} equations, it needs to be reformulated in this context. 

The first statement of the Jensen-Ishii's Lemma is due to Ishii \cite{ishii89} and relies on ideas developed by Jensen \cite{jensen88} when adapting the viscosity solutions approach to second-order elliptic equations. Let us recall that comparison principles themselves are, by many ways, the cornerstone of this theory since they provide not only uniqueness but also existence of continuous solutions when coupled with the Perron's method adapted by Ishii \cite{ishii89} to the framework of fully nonlinear, possibly degenerate, elliptic equations. We refer to the ``Users' guide''  \cite{cil92} for an introduction and a general presentation of the whole viscosity solutions theory.

Motivated by applications to finance but also by an increasing number of other ones 
(physical sciences, mechanics, biological models \textit{etc.}\,), the theory has 
been almost immediately extended (1986) to the context of partial integro-differential 
equations (PIDE for short), \textit{i.e.} partial differential equations involving 
nonlocal operators 
such as L\'evy ones
\begin{equation}\label{op:levy}
\mathcal{I}_L [u](x) = \int_{\R^d}\, (u (x+z) - u (x) - \nabla u (x) \cdot z \un_B (z) ) \mu (dz) 
\end{equation}
where $\mu$ is a singular measure and $B$ is some ball centered at $0$. To the best of our knowledge, the first paper devoted to this extension is the one by Soner \cite{soner86} 
in the context of stochastic control
of jump diffusion processes.  Since then, a huge literature has grown up and it would
be difficult (and irrelevant with respect to our goals) to cite all papers. Instead,
let us describe the difficulties that were successively 
overcome. Following Soner's work, 
a quite general class
of integro-differential equations nonlinear with respect to the nonlocal operators
were developed by Sayah \cite{sayah}. At that time, it was not possible to deal with equations
involving second-order derivatives of $u$ such as
$$
\lambda u + H (x,u,Du) - \Delta u - \mathcal{I}_L [u](x) =0\quad \mbox{in } \R^d.
$$
In the case of bounded measures, Alvarez and Tourin \cite{at96} 
obtained quite general results for parabolic equations. In \cite{bbp97,pham98} 
for instance, several comparison results
were obtained in special cases for singular measures. A  first attempt is
made by Jakobsen and Karlsen \cite{jk06} to give general results applicable to 
second-order elliptic equations. 
In order to get an analogous of Jensen-Ishii's lemma,
the authors have to assume that solutions are subquadratic. This assumption is not always
relevant since, roughly speaking, the restriction of the behaviour of solutions at infinity
is related to the integrability of the singular measure away from the origin. For instance,
in \cite{bbp97}, solutions with arbitrary polynomial growth are considered and even for 
a system of PIDE.\medskip

To sum it up,  the difficulties involved by elliptic nonlinear PIDE are
\begin{itemize}
\item
the coupling of second-order derivatives and nonlocal terms,
\item
the singularity of the measure appearing in the nonlocal operator,
\item
the behaviour at infinity of solutions. 
\end{itemize}
The third difficulty is studied in details in \cite{ai06} in the special case of
semilinear parabolic PIDE  and we refer 
the reader to it for a detailed discussion.

The present paper is focused on the first two difficulties we listed above and so let us be a bit more specific about them. To do so, we recall that when proving a comparison principle for a standard local equation, the definition of viscosity solutions with test-functions has to be completed with the equivalent definition in terms of so-called sub and superjets (see next section for a definition); and the classical Jensen-Ishii's Lemma allows to build elements of limiting semi-jets which play the role of first and second-order derivatives at the maximum point (after the doubling of variables).
\medskip

But, on one hand, this method relies on the fact that one can pass to the limit in the equation in order to write it for this larger class of generalized derivatives and, on the other hand, such a technic cannot be directly applied in the context of PIDE since test-functions are used not only to give a weak sense to the first and second derivatives of the solution but also to give a sense to the nonlocal operator around the singularity. 
\medskip

In the present work, we try to present a general framework to deal with
degenerate elliptic integro-differential equations. In particular, we 
give different equivalent definitions of solutions, we state and prove
a general stability result and we propose a nonlocal version of Jensen-Ishii's
lemma for these equations in order to prove comparison principles. We provide a comparison result which will be for us the occasion to discuss the assumptions one has to impose on the nonlocal term(s).
\medskip

Finding a proper notion of viscosity solution for degenerate elliptic integro-differential
equation is an important problem, especially when second-order terms are involved. 
See for instance the recent paper of Arisawa \cite{arisawa06a}.
The first two equivalent definitions of viscosity solutions (Definitions~\ref{def:phi}
and \ref{def:phi-u}) we give are quite classical, even if we present them in an original
general framework in order to deal with solutions with arbitrary growth at infinity
(not only polynomial growth). Freely speaking, the first one consists in replacing
the solution by the test-function on the whole space while the second one consists
in replacing it only around the singularity of the measure in the nonlocal 
operator. A third definition (Definition~\ref{defi:limiting}) is
given and it is new. It combines the use of semi-jets and test-functions.
We also prove on a special case (even if such a result holds true in a general 
setting) that one can in fact, in some way, use the function $u$ on the whole
space and thus obtain a definition that only relies on semi-jets. 
A whole section is devoted to examples in order to illustrate and justify 
the general framework we introduce. 

We next explain how to pass to the limit in PIDE. In the viscosity solution context,
the proper limits are the half-relaxed ones. As remarked in \cite{bs02} in a sublinear
setting, dealing with nonlocal operators involve specific technical difficulties. 
We explain how to overcome them without restricted the behaviour at infinity of 
solutions. 

The nonlocal version of Jensen-Ishii's lemma
relies on an adapted inf-convolution procedure: the slope of the test-function
is taken into account and the infimum is localized (in order to deal with functions
with arbitrary growth). Since the statement of the lemma is quite technical,
we immediately derive two corollaries that can be used in most examples of 
comparison principles to derive proper viscosity inequalities. 

Keeping in mind our illustrative purpose, we state and prove a comparison
principle under quite general assumptions focusing our attention on the coupling
between the nonlocal term and the $x$-dependence of the Hamiltonian.
\medskip

In a forthcoming work, we would like to continue our investigation by
studying in details nonlocal operators on bounded domains, in order to clarify
for instance what are the equivalent conditions of Dirichlet and Neumann conditions
for nonlocal operators (several interpretations are already given in the literature). We refer to Arisawa \cite{arisawa06a, arisawa06b} for results in this direction.
\medskip

The paper is organized as follows. In Section~\ref{sec:visc}, 
we recall two equivalent definitions of viscosity solutions for PIDE
and we show how to combine  semi-jets and test-functions in order to 
get a third equivalent definition. In Section~\ref{sec:ex}, we give
examples of singular measures, nonlocal operators and PIDE.
The general stability result is presented in Section~\ref{sec:stab}.
In Section~\ref{sec:lem}, we state our nonlocal version of Jensen-Ishii's lemma. 
In Section~\ref{sec:comparison}, we apply this lemma to proving a model
comparison principle.

\medskip
\noindent \textbf{Acknowledments.}
The authors wish to thank Roland Seydel for his interesting comments 
on the first version of this article which leads to this (a priori better) corrected
version.

\medskip
\noindent \textbf{Notation.} 
The scalar product of $x,y \in \R^d$ is denoted by  $x\cdot y$ and the
Euclidian norm of $x$ is denoted by  $|x|$. 
The unit ball of $\R^N$ (with $N=d$ most of the time) is denoted by  $B$.
A ball of radius $r$ centered at the origin is denoted by  $B_r$.  
The Hessian matrix of a twice differentiable function $u$ is denoted by  $D^2 u$. 
The $N \times N$ (real) identity matrix is simply denoted by  $I$. 
The space of $N \times N$ symmetric matrices with real entries is denoted by  $\S$.

\section{Nonlocal operators and viscosity solutions of PIDE}
\label{sec:visc}

In order to emphasize the common features of the (monotone) PIDE to
which viscosity solution theory applies in a natural way, we are going 
to consider the case of general equations written under the form
\begin{equation}\label{model}
F\left(x,u,\nabla u,D^2 u, \mathcal{I}  [x,u]  \right) = 0 \quad\hbox{in  }\R^d \; ,
\end{equation}
where $F$ is a continuous function satisfying the local and nonlocal degenerate ellipticity
conditions (E) (see below).

%-----------
Unfortunately, this simple, general model equation does not cover all the interesting cases~: in particular, the cases of the Bellman equation arising in stochastic control  (cf. (\ref{bellm}) below) or the system studied in \cite{bbp97} cannot be written in this way. But the ideas described in the present paper can be extended and used readily in this more complex framework.
A more general model equation could be the following one
\begin{equation}\label{model-gen}
F\left(x,u,\nabla u,D^2 u,  
\{\mathcal{I}_\alpha [x,u]\}_{\alpha \in A} \right) = 0 \quad\hbox{in  }\R^d \; ,
\end{equation}
where $F$ is continuous and  $\{\mathcal{I}[x,u]\}_{\alpha \in A}$ is a family of nonlocal terms.
%---------

\subsection{Assumptions}\label{sec:fram}

We first recall the classical definitions of semicontinuous envelopes and 
half-relaxed limits. For a locally bounded function $u$, its 
lower semicontinuous (lsc for short) envelope $u_*$ and its 
upper semicontinuous (usc for short) one $u^*$ are defined as follows
$$
u_* (x) = \liminf_{y\to x} u (y) \quad , \quad u^* (x) = 
\limsup_{y \to x} u (y)\; .
$$
For a sequence $(\ze)_\eps $ of uniformly locally bounded functions in some space $\R^m$
$$
\limiinf \ze (x) = \liminf_{\hrl} \ze (y) \quad , \quad \limssup \ze (x) = 
\limsup_{\hrl} \ze (y)\; .
$$

In order to be more specific on the assumptions we are going to use for the nonlocal term, we first introduce a space of functions $\mathcal{C}$ which is, typically, a set of functions satisfying a suitable growth condition. We use the following type of assumptions.
\medskip

\noindent \textsc{Assumption (C).} Given an upper semicontinuous
function $R:\R^d \to \R$, $\mathcal{C}$ is the space of functions
$u$ such that there exists a constant $\bar{c}>0$ such that 
for all $x\in \R^d$
$$|u(x)|\le {\bar c}(1+R(x)).$$
\ 

We remark that functions $u$ of $\mathcal{C}$ are locally bounded and the maximum or the minimum of two functions 
of $\mathcal{C}$ lies in $\mathcal{C}$. 
Moreover, if $K \subset \R^d$ is a compact set and if $\phi \in C^2 (K)$, there exists a function 
$\psi \in \mathcal{C}\cap C^2(K)$ such that $\psi = \phi$ on the interior of $K$. All these properties 
will be used throughout the paper.

%---------
\begin{ex}
The simplest example of set $\mathcal{C}$ is the space of bounded functions.
Another example is the space of sublinear functions, \textit{i.e.} functions $u$ for
which there exists a constant ${\bar c}$ such that
$$ 
|u (x)|\le {\bar c}( 1 + |x|). 
$$ 
We will see below that this kind of growth conditions is mainly related to the behaviour at infinity 
of the measure $\mu$ appearing in the nonlocal term. If this measure has a compact support, then we can even deal with $\mathcal{C} = C(\R^d)$.
\end{ex}
%---------

\medskip Our assumptions on the nonlocal term are the following.

\noindent \textsc{Assumption (NLT).}  For any $\delta >0$, there exist operators 
$\mathcal{I}^{1,\delta}[x,\phi]$, $\mathcal{I}^{2,\delta}[x,p,\phi]$ which are well-defined for 
any $x \in \R^d$ and $\phi \in \mathcal{C} \cap C^2 (\R^d)$ and which satisfies
  \begin{itemize}
  \item For any $x \in \R^d$ and $\phi \in \mathcal{C} \cap C^2 (\R^d)$,
 $ \mathcal{I}[x,\phi ] = \mathcal{I}^{1,\delta}[x,\phi] +
  \mathcal{I}^{2,\delta}[x,\nabla \phi (x),\phi]$. Moreover, for any $a \in \R$, $\mathcal{I}^{1,\delta}[x,\phi + a] =\mathcal{I}^{1,\delta}[x,\phi]$ and  $\mathcal{I}^{2,\delta}[x,\nabla \phi (x),\phi + a]=\mathcal{I}^{2,\delta}[x,\nabla \phi (x),\phi]$.
\item There exists $R_\delta>0$ with $R_\delta \to 0$ as $\delta \to 0$, such that if $\phi_1=\phi_2$ on $B(x,R_\delta)$
(resp. on $\R^d \setminus B(x,R_\delta)$), then
 $\mathcal{I}^{1,\delta}[x,\phi_1]=\mathcal{I}^{1,\delta}[x,\phi_2]$
(resp.  $\mathcal{I}^{2,\delta}[x,p,\phi_1]=\mathcal{I}^{2,\delta}[x,p,\phi_2]$).
  \item For any $\phi \in C^2 (\R^d)$ and $u \in \mathcal{C}$ 
such that $u - \phi$ attains a  maximum at $x$ on $B(x,R_\delta)$, 
there exists $\phi_k \in \mathcal{C}\cap C^2 (\R^d)$ such that:
\begin{eqnarray*}
u - \phi_k \mbox{ attains a global maximum at } x, \\
\begin{array}{l}
\mathcal{I}^{1,\delta}[x,\phi_k]\to\mathcal{I}^{1,\delta}[x,\phi],\\
\mathcal{I}^{2,\delta}[x,\nabla \phi_k (x), \phi_k]\to\mathcal{I}^{2,\delta}[x,\nabla \phi(x),u]
\end{array}\mbox{~as }k \to + \infty.
\end{eqnarray*}
  \item The operator $\mathcal{I}^{1,\delta}[x,\phi]$ is well-defined for any $x \in \R^d$ 
and $\phi \in C^2 (B(x,r)) \cap \mathcal{C}$ for any $r \in (0, R_\delta)$; moreover $\mathcal{I}^{1,\delta}[x,\phi] \to 0$ when $\delta \to 0$ and $\mathcal{I}^{1,\delta}[x_k,\phi_k] \to \mathcal{I}^{1,\delta}[x,\phi]$
 if $x_k \to x$ and $\phi_k \to \phi$ in $C^2 (B(x,r)) \cap C(\overline{B(x,R_\delta)}) $.
  \item The operator $\mathcal{I}^{2,\delta}[x,p,\phi]$ is defined for any $x \in \R^d$ and 
$\phi \in \mathcal{C}$. Moreover, if $x_k \to x$, $p_k \to p$ and  $(\phi_k)_k$ is a sequence
of uniformly locally bounded functions such that $|\phi_k| \le \psi$ 
 with $\psi \in \mathcal{C}$, 
$$\limsup_{k\to +\infty} 
\mathcal{I}^{2,\delta}[x_k,p_k,\phi_k] \le 
\mathcal{I}^{2,\delta}[x,p,\overline \phi] \quad (\hbox{resp.  }
\liminf_{k\to +\infty} 
\mathcal{I}^{2,\delta}[x_k,p_k,\phi_k] \ge \mathcal{I}^{2,\delta}[x,p,\underline \phi])
$$ 
where $\overline \phi := \limssup \phi_k$ 
(resp. $\underline \phi := \limiinf \phi_k$). 
\end{itemize}
%-----------
\begin{rem}
In the general case, \textit{i.e.} as far as \eqref{model-gen} is concerned,
Assumption~(NLT) must be satisfied by all the nonlocal operators 
$\mathcal{I}_\alpha$ with the same space $\mathcal{C}$. Such an
assumption is natural since only one singular measure appears in the
different nonlocal operators. 
\end{rem}
\begin{ex}
Consider the L\'evy operator appearing in \eqref{op:levy} with 
$\mu (dz) = \frac{dz}{|z|^{N+\alpha}}$. In this case for any $\delta >0$, 
\begin{eqnarray*}
\mathcal{I}^{1,\delta}[x, \phi]
&=& \int_{|z|\le \delta}(\phi(x+z) - \phi (x) - \nabla 
\phi (x) \cdot z \un_B(z)) \mu (dz),  \\
&=& \int_{|z|\le \delta}(\phi(x+z) - \phi (x) - \nabla 
\phi (x) \cdot z ) \mu (dz),  \\
\mathcal{I}^{2,\delta}[x,p,\phi]&=& 
 \int_{|z|\ge \delta}(\phi(x+z) - \phi (x)-p\cdot z \un_B (z) ) \mu (dz)\\
&=& \int_{|z|\ge \delta}(\phi(x+z) - \phi (x)) \mu (dz)
\end{eqnarray*}
(in addition we have used here the fact that $\mu$ is odd). If $\alpha >1$, sublinear functions are
integrable away from the origin and one can look for solutions with sublinear
growth at infinity. 
\end{ex}
%----------

Finally, we assume that $F$ is a continuous function satisfying the 
\emph{ellipticity assumption}:
\medskip

\noindent \textsc{Assumption (E).} For any $x \in \R^d$, $u \in \R$, 
$p \in \R^d$, $M, N \in \S$, $l_1, l_2 \in \R$
$$ F(x,u,p,M,l_1) \leq F(x,u,p,N,l_2)\quad\hbox{if } M\geq N,\ l_1 \geq l_2\; .$$
As we point it out making such an assumption, the fact that $F(x,u,p,M,l)$ 
is nonincreasing in $l$ is indeed part of the  ellipticity assumption on $F$.

\begin{rem}
As far as \eqref{model-gen} is concerned, $F$ is assumed to be
nondecreasing with respect to all nonlocal operators. 
\end{rem}

\subsection{Different definitions for viscosity solutions of PIDE}

In this subsection, we follow \cite{sayah}~by giving 
several definitions of viscosity solutions for Eq.~\eqref{model}.  
We will next prove that they are equivalent. Let us now give a first
definition of viscosity solution for \eqref{model}.  
%-----------------
\begin{defi}[Viscosity sub and supersolutions]\label{def:phi}
An usc function $u \in \mathcal{C}$ is a \emph{viscosity subsolution} of \eqref{model}
if, for any test-function $\phi \in \mathcal{C} \cap C^2 (\R^d)$, 
if $x$ is a global maximum point of $u - \phi$, then
$$
F(x,u(x), \nabla \phi (x), D^2 \phi (x), \mathcal{I}[x,\phi] ) \le 0 \; .
$$  
A lsc semicontinuous function $v \in \mathcal{C}$ is a \emph{viscosity supersolution} of \eqref{model}
if, for any test-function $\phi \in \mathcal{C} \cap C^2 (\R^d)$, if $x$ 
is a global minimum point of $u - \phi$, then
$$
F(x,u(x), \nabla \phi (x), D^2 \phi (x), \mathcal{I}[x,\phi] ) \ge 0 \; .
$$  
\end{defi}
%-------------------
\begin{rems}
\begin{itemize}
\item
It is also worth pointing out that, in Definition~\ref{def:phi}, we can 
as well assume that $\phi$ is $C^2$ in a small neighborhood of $0$ and is only continuous 
outside this neighborhood. This a consequence of the fourth and fifth points in Assumption (NLT).
\item 
One can develop a theory with subsolutions (resp. supersolutions) merely
 bounded from above (resp. from below) by a function of $\mathcal{C}$. 
Such an idea is somehow used in Proposition~\ref{prop:u-everywhere}.
\end{itemize}
\end{rems}
In the examples we present below, it will be clear (at least we hope so) 
that nonlocal terms $\mathcal{I}[x,w]$ are, in general, only defined  
for smooth functions $w$ because, typically, of the singularity of the L\'evy 
measure at $0$ and for functions with a suitable growth at infinity; this definition takes care 
of these two difficulties by using a test-function $\phi \in \mathcal{C} \cap C^2 (\R^d)$
 where we recall that $\mathcal{C}$ encodes the growth information.
\medskip

Assumption (NLT) is made in order that the following definition is equivalent to the 
previous one. 
%----------------
\begin{defi}[Viscosity sub and supersolutions]\label{def:phi-u}
An usc function $u \in \mathcal{C}$ is a \emph{viscosity subsolution} of \eqref{model}
iff, for any test-function $\phi \in C^2 (\R^d)$, if $x$ is a 
 maximum of $u - \phi$ on $B(x,R_\delta)$, then
$$
F(x,u(x), \nabla \phi (x), D^2 \phi (x), \mathcal{I}^{1,\delta}[x,\phi]
+ \mathcal{I}^{2,\delta}[x,\nabla \phi (x),u]  ) \le 0 \; .
$$  
A lsc function $v \in \mathcal{C}$ is a \emph{viscosity supersolution} of \eqref{model}
iff, for any test-function $\phi \in C^2 (\R^d)$, 
if $x$ is a minimum of $u - \phi$ on $B(x,R_\delta)$, then
$$
F(x,u(x), \nabla \phi (x), D^2 \phi (x), \mathcal{I}^{1,\delta}[x,\phi]
+ \mathcal{I}^{2,\delta}[x,\nabla \phi(x),v]) \ge 0 \; .
$$  
\end{defi}
%------------------
We now turn to a third, less classical definition, where we mix test-functions and 
sub-superjets. We first recall the definition of sub and superjets.
%-----------------------
\begin{defi}[Subjets and superjets] Let $u: \R^d \to \R$ be an usc function and $v: \R^d \to \R$ 
be a lsc function.\\
(i) A couple $(p,X) \in \R^d \times \S$ is a \emph{superjet} of $u$ at $x \in \R^d$  if
$$
u (x+z) \le u (x) + p \cdot z + \frac12 X z \cdot z + o(|z|^2). 
$$ 
(ii) A couple $(p,X) \in \R^d \times \S$ is a \emph{subjet} of $v$ at $x \in \R^d$ if
$$
u (x+z) \ge u (x) + p \cdot z + \frac12 X z \cdot z + o(|z|^2). $$
(iii) A couple $(p,X) \in \R^d \times \S$ is a \emph{limiting superjet} of $u$ 
at $x$ if there exists $(x_n,p_n,X_n) \to (x,p,X)$ such that $(p_n,X_n)$
is a superjet of $u$ at $x_n$ and $u(x_n) \to u(x)$. \\
(iv) A couple $(p,X)\in \R^d \times \S $ is a \emph{limiting subjet} of $v$ 
at $x$ if there exists $(x_n,p_n,X_n) \to (x,p,X)$ such that $(p_n,X_n)$
is a subjet of $v$ at $x_n$ and $u(x_n) \to u(x)$.
\end{defi} 
%------------------------
We will denote respectively by $J^+ u (x)$, $\overline{J}^+ u (x)$ the 
set of superjets and limiting superjets of $u$ at $x$ and by 
$J^- v (x)$, $\overline{J}^- v (x)$ the set of subjets and limiting subjets of $u$ at $x$.
\medskip

%----------------------------------------------------------------------
\begin{defi}[Sub-supersolutions and super-subjets] \label{defi:limiting}
A usc function $u \in \mathcal{C}$ is a \emph{viscosity subsolution} of \eqref{model}
if, for any test-function $\phi \in  C^2 (\R^d)$, if $x$ is a 
maximum point of $u - \phi$ on $B(x,R_\delta)$ and if $(p,X) \in J^+ u (x)$ with $p=D\phi(x)$ 
and $X \le D^2 \phi (x)$, then
$$
F(x,u(x), p, X, \mathcal{I}^{1,\delta}[x,\phi] 
+ \mathcal{I}^{2,\delta}[x,\nabla\phi(x),u] ) \le 0 \; .
$$  
A lsc function $v \in \mathcal{C}$ is a \emph{viscosity supersolution} of \eqref{model}
if, for any test-function $\phi \in C^2 (\R^d)$, 
if $x$ is a minimum point of $u - \phi$ on $B(x,R_\delta)$ and if $(q,Y) \in J^- v (x)$ 
with $q=D\phi(x)$ and $Y \ge D^2 \phi (x)$, then
$$
F(x,u(x), q, Y , \mathcal{I}^{1,\delta}[x,\phi] 
+ \mathcal{I}^{2,\delta}[x,\nabla\phi(x),v]) 
\ge 0 \; .
$$  
\end{defi}
%----------------------------------------------------------------------------
At first glance, this definition seems useless or at least rather strange; 
as Section~\ref{subsec:appl-lem} will show it, this is however the type of situation we face 
when applying the non local Jensen-Ishii's Lemma of Section~\ref{sec:lem}.
%-----------
\begin{prop}\label{prop:equivdef}
Definitions~\ref{def:phi},\ref{def:phi-u} and \ref{defi:limiting} are equivalent. 
\end{prop}
%---------
\begin{proof} 
We already justified that the first two definitions are equivalent. Let us prove
that so are Definitions~\ref{def:phi-u} and \ref{defi:limiting}.
We do it only for the subsolution case, the supersolution one being 
treated analogously. Changing $\phi$ in $\phi + \chi$ where $\chi$ is a $C^\infty$,
 positive function with compact support and such that $\nabla \chi (x) =0$, 
$D^2 \chi (x) = \alpha I$ with $\alpha >0$, we can assume that 
$X \leq D^2 \phi (x) - \alpha I$. Translating also $\phi$, 
we can also assume that $\phi(x) = u(x)$.

By classical results, there exists a smooth function $\psi: \R^d \to \R$ 
such that $\psi (x) = u(x)$,  $\nabla \psi (x) = p$, $D^2\psi (x) = X$ and 
$\psi \geq u$ in $\R^d$. We deduce from these properties that $x$ is 
a maximum of $u - \min(\psi,\phi)$ on $B(x,R_\delta)$ but, since $\psi (x) = \phi (x)$,
  $\nabla \psi (x) = \nabla \phi (x) $,  $D^2\psi (x) \leq D^2 \phi (x) - \alpha I$, 
we are sure that $\min(\psi,\phi) = \psi$ in a neighborhood of $0$. By 
Assumption~(NLT) (in particular 
the fact that the test-function needs to be smooth only in a neighborhood of $0$), 
we can use Definition~\ref{def:phi} to obtain
$$
F(x,u(x), \nabla \psi (x), D^2 \psi (x) , \mathcal{I}^{1,\delta}[x,\min(\psi, \phi) ] 
+ \mathcal{I}^{2,\delta}[x,\nabla\phi(x),u] ) \le 0 \; .
$$  
By using the monotonicity of $F$ in $l$, we can as well replace 
$\min(\psi, \phi)$ by $\phi$, and the proof is complete.
\end{proof}

\subsection{An additional proposition}

In \cite{imbert05}, general test-functions are replaced with functions
$ 
\phi (t,x) = \alpha t + p \cdot x \pm \sigma |x|^2 + o (|t|+|x|^2).
$
In particular, it is noticed that one can choose $\delta=0$ 
in Definition~\ref{def:phi-u} for such test-functions.  
Let us explain this in a general setting. 

For simplicity, consider a bounded viscosity subsolution $u$ of a special case of \eqref{model}
\begin{equation}\label{model-2}
F (x,u, \nabla u, D^2 u, \mathcal{I}[u]) =0 \quad \mbox{~in~}\R^d
\end{equation} 
where, if $\mu$ denotes a singular odd measure,
$$
\mathcal{I}[u](x) = \int (u(x+z) - u (x)- \nabla u (x)\cdot z \un_B(z))\mu (dz).
$$
For such a nonlocal operator, we choose $R_\delta=\delta \le 1$ and 
\begin{eqnarray*}
\mathcal{I}^{1,\delta}[x,\phi]=  \int_{|z|\le \delta} (\phi(x+z) - \phi (x)
- \nabla \phi (x)\cdot z )\mu (dz),\\
\mathcal{I}^{2,\delta}[x,p,u]=  \int_{|z|\ge \delta} (u(x+z) - u (x)
- p\cdot z \un_B(z))\mu (dz).
\end{eqnarray*}
If $\phi \in C^2 (\R^d)$ is a test-function such that $u-\phi$ attains a global
maximum at $x$ and $p=\nabla \phi (x)$, we have
$$ \forall z \in B_\delta, \quad u (x+z) - u (x) -p \cdot z \le \phi (x+z) -\phi (x) - p \cdot z \; ,$$
and therefore $\psi (z) :=[ u (x+z) - u (x) -p \cdot z] - [ \phi (x+z) -\phi (x) - p \cdot z] $ is nonpositive. This implies that $ \int_{B_\delta} \psi(z) \mu (dz)$ is well-defined and since $z \mapsto \phi (x+z) -\phi (x) - p \cdot z
\in L^1 (B_\delta, d\mu)$, the integral $ \int_{B_\delta} (u (x+z) - u (x) -p \cdot z) \mu (dz)$ is well-defined too. Moreover, by the monotone convergence theorem, we can pass to the limit in the equality
$$\int_{r \le |z|\le \delta} (u (x+z) - u (x) -p \cdot z) \mu (dz) = \int_{r \le |z|\le \delta}( \phi(x+z) - \phi(x)  
- p \cdot z) \mu (dz) -  \int_{r \le |z|\le \delta} \psi(z) \mu (dz)\;
$$
and, consequently, one can define 
$$\mathcal{I}^{2,0} [x,p,u] = \lim_{\delta \to 0+} \mathcal{I}^{2,\delta} [x,p,u] 
\in \{-\infty \}\cup \R.$$ 
The (nonincreasing) function $F$ has a limit $l \in \R \cup \{ + \infty\}$ as
goes to $-\infty$. But since
\begin{equation} \label{eq:obtenue}
F(x,u(x), \nabla \phi (x), D^2 \phi (x), \mathcal{I}^{1,r} [x,\phi]  + \mathcal{I}^{2,\delta}
[x,p,u]) \le 0,
\end{equation}
we conclude that
$$
F(x,u(x), \nabla \phi (x), D^2 \phi (x), \mathcal{I}^{2,0} [x,p,u]) \le 0.
$$
Indeed, if $F (\lambda) \to + \infty$ as $\lambda \to -\infty$, then \eqref{eq:obtenue} implies
that $\mathcal{I}^{2,0} [x,p,u]$ is finite and passing to the limit as $\delta\to 0$ gives
the result. And if $F(\lambda) \to l_{-\infty} \in \R$ as $\lambda \to -\infty$, 
$F(x,u(x), \nabla \phi (x), D^2 \phi (x), l) \le 0$ for any $l \in \R\cup \{-\infty\}$. 
It is also clear that such a discussion can be adapted to the case of supersolutions. 
Let us sum up this discussion in the following proposition. 
%----------------
\begin{prop} \label{prop:u-everywhere}
For any subsolution $u$ of \eqref{model} and any test-function $\phi \in C^2 (\R^d)$
such that $u -\phi$ attains a global maximum at $x$,  
\begin{eqnarray*}
\mathcal{I}^{2,0} [x,\nabla \phi (x),u] \in  \{-\infty \}\cup \R \\
F(x,u(x), \nabla \phi (x), D^2 \phi (x), \mathcal{I}^{2,0} [x,\nabla\phi(x),u] ) \le 0.
\end{eqnarray*}
Moreover, if $F \to + \infty$ as $l \to -\infty$, then 
$\mathcal{I}^{2,0} [x,\nabla\phi (x),u]$ is finite.  
\end{prop}
%---------------
\begin{rems}
\begin{itemize}
\item
For clarity, we only treated a special case but the attentive reader can check
that such a proposition holds true for all the examples we will give below. 
\item
Remark that one can use this proposition to give an alternative
proof of the fact that Definition~\ref{defi:limiting} is equivalent
to Definitions~\ref{def:phi}~and \ref{def:phi-u}. 
\end{itemize}
\end{rems}
Following this idea, Arisawa \cite{arisawa06a} 
considered quadratic test-functions
$$
\alpha t + p \cdot x + \frac12 A x \cdot  x+ o (|t|+ |x|^2)
$$
 and keep $r>0$. Since this definition
involves some technicalities when proving comparison principles, we will not use this
fourth definition.

\section{Examples}
\label{sec:ex}

PIDE's involve nonlocal operators that are (partially) characterized
by an (eventually singular) measure. In this section, we give examples
of such measures, operators and equations. Let us start with measures. 

\subsection{Singular measures}

In applications, typically in stochastic control with jump processes,
positive singular measures appear in nonlocal operators. For instance, the pure
jump part of a L\'evy process is characterized by such a measure that is
refered  to as the L\'evy measure. It consequently appears in the associated 
infinitesimal operator one has to deal with in the associated Bellman equation. 

Let us next give two examples of L\'evy measures
\begin{eqnarray}\label{mes:alpha}
\mu_1 (dz) = g\left( \frac{z}{|z|} \right) \frac1{|z|^{N+\alpha}} dz  
\quad \mbox{ with }\alpha
\in (0,2),\;\hbox{ in  }\R^d,\\
\label{mes:gamma}
\mu_2 (dz) = ( \un_{(0,+\infty)} (z) e^{-\gamma^+ z}
+ \un_{(-\infty,0)}(z)e^{\gamma^- z})\frac1{|z|} dz,\;\hbox{ in  }\R.
\end{eqnarray}
The measure $\mu_1$ (resp. $\mu_2$) is an anisotropic
(resp. tempered) $\alpha$-stable L\'evy measure on $\R^d$ (resp. on $\R$).

We would like next to discuss how singular the measure is around the origin
and how it decreases at infinity. If one considers the measures $\mu_1$ and $\mu_2$ 
introduced previously
\begin{eqnarray*}
\int_B |z|^{\alpha'} \mu_1 (dz) <+\infty \quad \mbox{ for }\alpha' > \alpha,  \\
\int_{-1}^1 |z|^{\alpha'} \mu_2 (dz) <+\infty \quad \mbox{ for }\alpha' >0,\\
\int_{\R^d \setminus B} |z|^{\alpha'}\mu_1 (dz) < +\infty \quad \mbox{ for~}\alpha' < \alpha, \\
\int_{\R \setminus [-1,1]}|z|^p \mu_2 (dz) < +\infty \quad \mbox{for any } p \in \N.
\end{eqnarray*}

Let us mention that a L\'evy measure $\mu$ always satisfy
\begin{equation}\label{cond:levymeasure}
\int \min ( |z|^2, 1)\;  \mu (dz) < + \infty.
\end{equation}

\subsection{Nonlocal operators}

Now we examine typical examples of nonlocal operators appearing in the applications.
For instance, if  $\varphi$ denotes 
a smooth function satisfying suitable growth conditions, the following two nonlocal operators
appear in \cite{bbp97}
\begin{eqnarray}
\label{K}
K \varphi (x) = \int_{\R^d} (\varphi (x + \beta (x, z)) - \varphi(x) -  
\nabla \varphi (x) \cdot \beta (x, z) ) \lambda (dz) \\
\label{B}
B \varphi  (x) = \int_{\R^d} (\varphi  (x + \beta  (x, z) ) - \varphi (x))
\gamma  (x, z) \lambda  (dz);
\end{eqnarray}
as a matter of fact, we redefine the measure $\lambda$ appearing in the
compensator at the origin by setting:
$\lambda (\{0 \})=0$ so that the domain of integration is the whole space
$\R^d$ and not $E:= \R^d \setminus \{ 0 \}$. 
The functions $\beta, \gamma$ have to satisfy suitables assumptions.

The natural assumptions coming from the probabilistic formulation are \eqref{cond:levymeasure} 
and
$$  | \beta (x, z) |, |\gamma (x, z)| \leq  K|z| \quad \hbox{if  }|z|\leq 1\, .$$
It is easy to check that, under these assumptions, the $\mathcal{I}^{1,\delta}[x,\phi]$ 
terms, namely
$$ 
\int_{\{|z|\leq \delta\}} (\phi (x + \beta (x, z)) - \phi (x) -  
\nabla \phi (x) \cdot  \beta (x, z)  )\lambda  (dz) 
\quad\hbox{or}\quad \int_{\{|z|\leq \delta\}} (\phi (x + \beta (x, z)) 
- \phi (x) )\gamma (x,z) \lambda  (dz) 
$$
are well-defined for any $x\in \R^d$ and any smooth function $\phi$, and they satisfy (NLT).

On the other hand, the $\mathcal{I}^{2,\delta}[x,p,\phi]$ terms, namely
$$ 
\int_{\{|z| \geq \delta\}} (\phi (x + \beta (x, z)) - \phi (x) -  
p \cdot  \beta (x, z)  )\lambda  (dz) 
\quad\hbox{or}\quad \int_{\{|z|\geq \delta\}} 
(\phi (x + \beta (x, z)) - \phi (x) )\gamma (x,z) \lambda  (dz) 
$$
are well-defined and satisfy (NLT) if $\beta, \gamma$ have suitable integrability 
properties w.r.t the measure $\lambda (dz)$; these integrability properties determine
 also the space $\mathcal{C}$ on which the operators $\mathcal{I}^{2,\delta}[x,p,\phi]$ 
are defined. In \cite{bbp97}, $\beta$ is assumed to be bounded and therefore (NLT) is 
readily satisfied for $\mathcal{C} = C(\R^d)$.
\medskip

In stochastic control with jump processes (see for instance \cite{os}), 
one can consider the so-called L\'evy-Ito
diffusions whose infinitesimal generators are of the form:
\begin{equation}
\label{op:levy-ito} \mathcal{I}_{LI} [u](x) = \int ( u (x+j (x,z) ) 
- u (x) - \nabla u (x) \cdot j(x,z) \un_B (z) ) \mu (dz)
\end{equation}
where $\un_B$ denotes the indicator function of the ball $B$,
$\mu$ is a L\'evy measure (hence it satisfies \eqref{cond:levymeasure}) and 
$j (x,z)$ is the size of the jumps at $x$. In order that the operator is well-defined, 
one assumes:  
\begin{equation} \label{cond:j}
|j(x,z) |\le {\bar c}|z|. 
\end{equation}
The simplest example of $j$ function is $j(x,z) = z$ and in this case, operators
are the infinitesimal generators of any pure jump L\'evy process; hence, it is referred
to as L\'evy operators (see \cite{bertoin} for details). 
One can check that (NLT) is also satisfied with   
$\mathcal{C} = C(\R^d)$ and:
\begin{eqnarray*}
\mathcal{I}^{1,\delta}_{LI}[x,\phi]=\int_{|z|\le \delta} 
( \phi(x+j(x,z))-\phi(x)-\nabla\phi(x)\cdot j(x,z)\un_B(z))\mu(dz), \\
\mathcal{I}^{2,\delta}_{LI}[x,p,\phi]=\int_{|z|\ge \delta} 
(\phi(x+j(x,z))-\phi(x)-p\cdot j(x,z)\un_B(z))\mu(dz).
\end{eqnarray*}
\medskip

Let us conclude this subsection by giving other examples one can find in the
(huge) literature concerning nonlocal operators:
\begin{eqnarray}
\label{op:silv}
\mathcal{I}_{Si}[x,u] (x) 
= \int ( u (x+z) - u(x) -\nabla u (x) \cdot z \un_B (z) ) ) K(x,z) dz \\ 
\label{op:sayah}
\mathcal{I}_{Sa} [x,u] (x) = \int ( u (x+z) - u(x) -\nabla u (x) \cdot z \un_B (z) ) ) 
\mu (x,dz) \\
\label{op:gm}
\mathcal{I} [x,u](x) = \int (u (x+j(x,z)) - u (x) - \nabla u (x) \cdot j (x,z) 
\un_B (z) ) \mu (x,dz) 
\end{eqnarray}
where $\mu(x,dz)$ are bounded measures that are possibly singular at the origin. 
Operators
appearing in \eqref{op:silv} are considered in \cite{silvestre06}; see this paper for further
details. Sayah \cite{sayah} considered operators of the
form~\eqref{op:sayah}. To finish with, operators of the form \eqref{op:gm} are probably 
the most general ones and they appear for instance in \cite{gm02}.

\subsection{Examples of PIDE}

A first example of (semilinear) PIDE is one appearing in the context of
growing interfaces \cite{woyczynski}
$$
\partial_t u + \frac12|\nabla u|^2 - \mathcal{I}_L [u] =0
$$
where $\mathcal{I}_L$ is a L\'evy operator. 

A second simple example of nonlinear PIDE 
is for instance the nonlinear diffusion arising in the context
of the homogenization for dislocation dynamics \cite{imr06}
$$
\partial_t u = \overline{H}(\nabla u, \mathcal{I}_L [u])
$$
where $\mathcal{I}_L$ is an anisotropic L\'evy operator of order $1$.

An important example of a nonlinear elliptic PIDE is the Bellmann equation
associated with a stochastic control problem
\begin{equation}\label{bellm}
\lambda u + \sup_{\alpha \in \mathcal{A}} \{ - \mathcal{I}^\alpha [u]- 
\frac12 D^2 u (x) \; \sigma_\alpha (x) \cdot \sigma_\alpha (x)  - b^\alpha (x) \cdot 
\nabla u (x) - f^\alpha (x) \} =0 
\end{equation}
where 
$$
\mathcal{I}^\alpha [u]=\int (u (x+j_\alpha (x,z) ) 
- u (x) - \nabla u (x) \cdot j_\alpha (x,z) \un_B (z) ) \mu (dz) .
$$

\section{Stability results}
\label{sec:stab}

In this section, we state and prove general stability results for viscosity solutions of PIDE. 
These results are a slight generalization of all previous analogous results \cite{sayah,at96,bs02} since it covers the case of general unbounded solutions.  To do so, we use the framework introduced in Section~\ref{sec:fram}.

For $\eps >0$, we consider sub or supersolutions $\ue$ of 
\begin{equation}\label{model-n}
\Fe (x,\ue, \nabla \ue, D^2 \ue, \mathcal{I}[x,\ue])=0 \quad \mbox{in }\R^d \; .
\end{equation}

The main stability result is the following one. 
\begin{thm}
Let $(\Fe)_\eps$ be a sequence of locally uniformly bounded nonlinearities satisfying the ellipticity condition (E) and let $(\ue)_\eps$ be a sequence of subsolutions (resp. supersolutions) of \eqref{model-n} such that there exists ${\bar c} >0$ such 
that
$$ |\ue (x)| \le {\bar c} (1+R(x)) \quad \mbox{in  }\R^d\; .$$
Then $\ou:= \limssup \ue$ (resp. $\uu:= \limiinf \ue$) is a subsolution 
(resp. supersolution) of \eqref{model} with $\uF:= \limiinf \Fe$ (resp. $\oF:= \limssup \Fe$).
\end{thm}

\medskip

The second stability result (which is partly a consequence of the first one) concerns the supremum/infimum of a family of sub/supersolutions. It is the cornerstone of Perron's method when proving the existence of a solution. 

To state it in a general way, we consider sub and supersolutions $u_\alpha$, for $\alpha \in \mathcal{A}$, of 
\begin{equation}\label{model-a}
\Fa (x,u_\alpha, \nabla u_\alpha, D^2 u_\alpha, \mathcal{I}[x,u_\alpha])=0 \quad \mbox{in }\R^d,
\end{equation}
where $\mathcal{A}$ can be any set.

\begin{thm}
Let $(F_\alpha)_{\alpha \in \mathcal{A}}$ be a family of uniformly locally bounded from
above (resp. from below) nonlinearities satisfying the ellipticity condition (E) and let 
$(u_\alpha)_{\alpha \in \mathcal{A}}$ be a family of  
subsolutions (resp. supersolutions) of \eqref{model-a} such that there exists ${\bar c} >0$ such 
that, for any $\alpha \in\mathcal{A}$ and $x\in \R^d$
\begin{equation}\label{cond:lim}
u_\alpha (x) \le {\bar c} (1+R(x)) \qquad \mbox{(resp. $u_\alpha (x) \ge -{\bar c} (1+R(x))$).}
\end{equation}
We set $u=\sup_{\alpha \in \mathcal{A}}u_\alpha$ (resp. $v=\inf_{\alpha \in \mathcal{A}}u_\alpha$).
Then $u^*$ (resp. $v_*$) is a subsolution (resp. supersolution) of \eqref{model}
where $\uF=\left(\inf_{\alpha \in \mathcal{A}} F_\alpha\right)_*$ 
(resp. $\oF=\left(\sup_{\alpha \in \mathcal{A}} F_\alpha\right)^*$). 
\end{thm}

\begin{rems}
\begin{itemize}
\item
The reader can check that both stability results can be easily adapted to the
general case of several nonlocal operators such as \eqref{model-gen}. 
\item
The condition~\eqref{cond:lim} is not necessary for local equations but
it cannot be avoided for nonlocal ones. A special case of it appears in 
\cite{bs02}. 
\end{itemize}
\end{rems}
The proofs of both theorems are easy adaptations of classical ones. 
For the sake of completeness, we give a proof of the first one and
we let the reader check that the classical proof for the first one can be also
adapted. 
\begin{proof}
We only provide the proof for $\ou$ since the other case follows along
the same lines. 
In order to prove that $\ou$ is a subsolution of the limit equation, we consider a test-function
$\phi$ and a maximum point $x$ of $\ou -\phi$ on $B(x,R_\delta)$ (see Assumption~(NLT)). 

First we can assume without loss of generality that $x$ is a strict maximum point of
 $\ou -\phi$ on $B(x,R_\delta)$: indeed we may replace $\phi$ by 
$\tilde\phi=\phi + \alpha \chi$ where $\chi$ is a $C^\infty$ function 
whose support is $B(x,2R_\delta)$ and such that $\chi >0$ on $B(x,R_\delta)\setminus \{x\}$, 
$\chi(x)=0$, $\nabla \chi (x) =0$, 
$D^2 \chi (x) =0$.

Next we consider a subsequence such that
$$
\ou (x) = \lim_{\eps'} u_{\eps'} \left(x_{\eps'}\right).
$$ 
Since $x$ is a strict maximum point of $\ou -\phi$ on $B(x,R_\delta)$, classical arguments show that $u_{\eps'} -\tilde\phi$ attains a maximum on $B(x,R_\delta)$ at $y_{\eps'} \in B(x,R_\delta)$ ; moreover 
$$
x = \lim_{\eps'} y_{\eps'} \quad \hbox{and}\quad \ou (x) = \lim u_{\eps'} (y_{\eps'}).
$$
Since $u_{\eps'}$ is a subsolution of \eqref{model-n}, we have
$$
F_{\eps'} (x_{\eps'}, u_{\eps'} (x_{\eps'}), \nabla \tilde\phi(x_{\eps'}),D^2\tilde\phi(x_{\eps'}),l_{\eps'})\le 0
$$
with 
$$
l_{\eps'}=\mathcal{I}^{1,\delta}[x_{\eps'},\tilde\phi]
+\mathcal{I}^{2,\delta}[x_{\eps'},\nabla\tilde\phi(x_{\eps'}),u_{\eps'}].
$$
By (NLT) and \eqref{cond:lim}, we can conclude that
$$
\limsup_{\eps'}\, l_{\eps'} \le \mathcal{I}^{1,\delta}[x,\tilde\phi]
+\mathcal{I}^{2,\delta}[x,\nabla\phi(x),\ou]
$$
and the definition of $\uF$ together with the nonlocal ellipticity permits to get
$$
\uF (x,\ou(x),\nabla\phi(x),D^2\phi(x),\mathcal{I}^{1,\delta}[x,\tilde\phi]
+\mathcal{I}^{2,\delta}[x,\nabla\phi(x),\ou])\le 0. 
$$ 
Finally we let $\alpha$ tend to $0$ and the proof is complete. 
\end{proof}

\section{A nonlocal version of Jensen-Ishii's Lemma}
\label{sec:lem}

In order to state our result, we first introduce the inf and sup-convolution 
operations we are going to use.

\subsection{Modified inf/sup-convolution procedures}
 
For any usc function $U: \R^m \to \R$ and any lsc function $V: \R^m \to \R$, we set
$$
R^\alpha [U] (z,r) = \sup_{|Z-z| \le 1} \bigg\{ U (Z) 
- r  \cdot (Z-z)  - \frac{|Z-z|^2}{2 \alpha} \bigg\} $$
$$
R_\alpha [V] (z,r) = \inf_{|Z-z| \le 1} \bigg\{ V (Z) 
+ r  \cdot (Z-z)  + \frac{|Z-z|^2}{2 \alpha} \bigg\} $$
Notice that $R_\alpha [V] = - R^\alpha [-V].$

%----------------
\begin{prop} \label{prop:infconv}
For any usc function $U: \R^m \to \R$ and any lsc function $V: \R^m \to \R$, the function
 $R^\alpha[U], R_\alpha [V] $ satisfy the following properties 
\begin{enumerate}
\item For any $x, r \in \R^m$, $R^\alpha[U] (x, r) \ge U$ and $R_\alpha[V] (x, r) \le V$
\item For any $x \in \R^m$ and $ {\bar k} >0$, there exists $\bar \alpha = \bar \alpha (x, {\bar k})$ such that,
 for $0<\alpha \leq \bar \alpha$, $R^\alpha[U](\cdot,r)$ is semi-convex in $B(x, {\bar k})$ 
(resp. $R_\alpha[V](\cdot,r)$ is semi-concave in $B(x, {\bar k}))$
\item Assume that $U \in C^2 (\R^m)$ (resp. $V \in C^2 (\R^m)$). For any $x \in \R^m$ and $ {\bar k} >0$,
 there exists $\bar \alpha = \bar \alpha (x, {\bar k})$ such that, for $0<\alpha \leq \bar \alpha$, 
then $R^\alpha [U]$   (resp. $R_\alpha[V] $) is $C^2$ in $B(0, {\bar k})$. Moreover, $R^\alpha [U]$ (resp.
$R_\alpha [V]$) converges towards $U$ (resp. $V$) in $C^2(B(0, {\bar k}))$ as $\alpha \to 0$. 
\item
If $\displaystyle R^\alpha [U](z,r) = U (\bar{z}) - r \cdot (\bar{z}-z) 
- \frac{|\bar{z}-z|^2}{2 \alpha}$ and if $|\bar{z}-z| < 1$,
\begin{eqnarray}
(s,A) \in J^+ R^\alpha [U] (z,r) 
\Rightarrow (s,A) \in J^+ U (\bar{z}) \mbox{ and }s = r - \frac{z-\bar{z}}{\alpha},
\label{eq:inclus}\\
(r,A) \in \overline{D}^{2,+} R^\alpha[U] (z,r) \Rightarrow (s,A) \in \overline{D}^{2,+} U (z).
\label{eq:alalim}
\end{eqnarray}
\end{enumerate} 
\end{prop}
%--------
\begin{proof}
The first two points are clear. The third point is a consequence of 
the analogous classical result: indeed, for $\alpha$ small enough, 
the supremum and infimum are achieved for $|Z-z|< 1$ and one can apply 
the same proof as for the classical sup and inf-convolutions 
(maximizing or minimizing w.r.t. $Z \in \R^m$ as if $U$ and $V$ were bounded). Let us focus on the fourth one.
 Eq.~\eqref{eq:inclus} is a simple adaptation  of the classical result about 
sup-convolution. Eq.~\eqref{eq:alalim} is a consequence
of it. Indeed, by definition of limiting superjets, there exists $(r_n,A_n) \in
J^+ R^\alpha [U] (z_n,r)$ such that $(r_n,A_n,z_n) \to (r,A,z)$. Moreover, by
\eqref{eq:inclus}, we have: $(r_n,A_n) \in J^+U (z_n + \alpha (r -r_n))$.
The proof is now complete. 
\end{proof}

\subsection{Statement and proof of the lemma}

We can now state our nonlocal version of Jensen-Ishii's lemma.
%--------
\begin{lem}[Nonlocal Jensen-Ishii's Lemma] \label{lem:global-jil}
Let $u$ and $v$ be respectively a usc and a lsc function defined on $\R^d$ and let $\phi$ be a 
$C^2$ function defined on $\R^{2d}$. If $(\xb,  \yb) \in \R^{2d}$ is a zero global maximum point
 of $u(x) -v(y) - \phi (x,y)$ and if $p:=D_x \phi (\xb,\yb)$, $q:=D_y \phi (\xb,\yb)$, then the
 following holds
 \begin{eqnarray}
\nonumber
u (x) - v (y) & \le & R^\alpha [u] (x,p) -R_\alpha [v] (y,-q)  \\
 & \le & R^\alpha [\phi]((x,y),(p,q)) \; ,\label{eq:interpos}\\
\label{eq:egal-u}
u(\xb) & = & R^\alpha [u] (\xb,p) , \\
\label{eq:egal-v}
v (\yb) & = & R_\alpha [v](\yb, -q) , \\
\label{eq:egal-phi}
R^\alpha [\phi] ((\xb,\yb),(p,q)) & = & \phi (\xb, \yb)\; .
\end{eqnarray}
Moreover, for any $ {\bar k}>0$, there exists $\bar \alpha ( {\bar k}) >0$ such that, for any 
$0 < \alpha \leq \bar \alpha ( {\bar k})$, we have: 
there exists sequences $x_k \to x$, $y_k \to y$, $p_k \to p$, $q_k \to q$,
 matrices $X_k,Y_k$ and a sequence of functions $\phi_k$, converging to the function 
$\phi_\alpha:=R^\alpha [\phi]((x,y),(p,q))$ uniformly in $\R^m$ and in $C^2 (B((\xb,\yb),  {\bar k}))$,
 such that
\begin{eqnarray}
(x_k,y_k)\hbox{ is a global maximum point of $u-v-\phi_k$}\\
u(x_k) \to u(\xb),\ v(y_k) \to v(\yb)\\
(p_k,X_k) \in J^+ u (x_k ) \\ 
(-q_k ,Y_k ) \in J^- v (y_k) \\
\label{eq:ineg-matrix}
- \frac1\alpha I \le \left[\begin{array}{ll} X_k & 0 \\ 0 & -Y_k \end{array} \right]
\le D^2\phi_k (x_k, y_k)\; .
\end{eqnarray}
Moreover $p_k=D_x \phi_k (x_k, y_k)$, $q_k=D_y \phi_k (x_k, y_k)$, and $\phi_\alpha (\xb,\yb) = \phi (\xb,\yb)$, $D \phi_\alpha (\xb,\yb)=D \phi (\xb,\yb)$.
\end{lem}
%--------
\begin{rem} The nonlocal Jensen-Ishii's Lemma is stated for
functions $u$ and $v$ which are defined in $\R^d$ but the same result
holds if $u$ and $v$ are defined only on a (closed) subset of $\R^d$.
Indeed, following the ``User's guide'' (cf. \cite{cil92}, p.~57), it
suffices to extend $u$ and $v$ in a suitable way outside this subset
(and typically $u$ by $-\infty$ and $v$ by $+\infty$). This remark is
important when one wants to deal with problems set in a domain of $
\R^d$ with boundary conditions.
\end{rem}
\begin{proof}
Eq.~\eqref{eq:interpos} is a simple consequence of Proposition~\ref{prop:infconv} 
of the fact that $u-v-\phi$ attains a zero global maximum. Eq.~\eqref{eq:egal-phi} 
is a consequence of the regularity of $\phi$ and more precisely of the property
$$ 
\phi (x,y) \leq \phi (\xb,\yb) + p\cdot (x-\xb) + q \cdot (y-\yb) + K |(x,y) - (\xb,\yb)|^2
$$
for some constant $K>0$ and for $|(x,y) - (\xb,\yb)| \leq 1$.
Eq.~\eqref{eq:egal-phi} implies \eqref{eq:egal-u} and \eqref{eq:egal-v}. 

The function $(x,y) \mapsto  R^\alpha [u](x,p) - R_\alpha [v] (y,-q) 
- R^\alpha [\phi ]((x,y), (p,q))$
 is semi-convex and achieves a global maximum point at $(\xb, \yb)$. 
In order to apply Lemma~A.3 of \cite{cil92}, we have to transform $(\xb, \yb)$ 
into a strict maximum point. To do so, we consider a smooth bounded function 
$\chi$ such that $\chi >0$ in $B \setminus\{(\xb,\yb)\}$
 and $\chi=0$ outside; $(\xb, \yb)$ is a strict 
maximum point of
$$
(x,y) \mapsto  R^\alpha [u](x,p) - R_\alpha [v] (y,-q) - 
R^\alpha [\phi ]((x,y), (p,q))-\delta \chi(x,y)\; ,
$$
for any $\delta >0$.

Next we consider a smooth function $\psi: \R^d \times \R^d \to \R$, 
with compact support and such that $\psi (x,y) = 1$ if $|(x,y) - (\xb,\yb)| \leq 1$. 
We are going to apply Lemma~A.3 of \cite{cil92} to the function
$$
(x,y) \mapsto  R^\alpha [u](x,p) - R_\alpha [v] (y,-q) 
- R^\alpha [\phi ]((x,y), (p,q))-\delta \chi(x,y) + \psi(x,y)\left (r\cdot x + s \cdot y\right)
\;,$$
for $r,s \in \R^d$ close to $0$. On one hand, the fact that $(\xb,\yb)$ 
is a strict global maximum of $R^\alpha [u](x,p) - R_\alpha [v] (y,-q) 
- R^\alpha [\phi ]((x,y), (p,q))-\delta \chi(x,y)$ implies that this 
function has global maximum points near $(\xb,\yb)$ for $r,s$ close enough to $0$, 
and, on the other  hand, since $\psi$ is $1$ in a neighborhood of $(\xb,\yb)$, 
we can readily apply Lemma~A.3 of \cite{cil92} in this neighborhood where the function is nothing but
$$
(x,y) \mapsto  R^\alpha [u](x,p) - R_\alpha [v] (y,-q) - 
R^\alpha [\phi ]((x,y), (p,q))-\delta \chi(x,y) +r\cdot x + s \cdot y\; .
$$
 Combining it with Theorem~A.2 of \cite{cil92}, we deduce that, for any $\delta >0$, 
there exist sequences $(r^\delta_m)_m$, $(s^\delta_m)_m$ converging to $0$ as $m \to +\infty$ and global maximum points
 $(x^\delta_m,y^\delta_m)_m$ of the above function such that $R^\alpha [u]$ is twice 
differentiable at $x^\delta_m$ and $R_\alpha [v]$ is twice differentiable at $y^\delta_m$. Of course, for fixed $\delta$, $(x^\delta_m,y^\delta_m)_m$
converges to $(\xb,\yb)$ as $m\to +\infty$ since $(\xb,\yb)$ is a strict global maximum of $R^\alpha [u](x,p) - R_\alpha [v] (y,-q) 
- R^\alpha [\phi ]((x,y), (p,q))-\delta \chi(x,y)$.

If now we choose $\delta = k^{-1}$, and we consider $(\tilde{x}_k, \tilde{y}_k) 
= (x^{1/k}_{m_k},y^{1/k}_{m_k})$ with $m_k$ chosen large enough so that 
$$
|(x^{1/k}_{m_k},y^{1/k}_{m_k}) - (\bar{x},\bar{y})| \le \frac1k 
$$
 and
$$ 
\tilde\phi_k (x,y) = R^\alpha [\phi ]((x,y), (p,q))+\frac1k \chi(x,y) 
-r^{1/k}_{m_k} \cdot x - s^{1/k}_{m_k} \cdot y\; ,
$$
then we have a sequence $(\tilde{x}_k, \tilde{y}_k)$ of maximum points of $R^\alpha [u](x,p) 
- R_\alpha [v] (y,-q) - \tilde \phi_k (x,y)$. One can remark that $\tilde \phi_k$ are small, 
smooth perturbations of $R^\alpha [\phi ]((x,y), (p,q))$. We also point out that we have chosen $m_k$ so that
$(\tilde{x}_k, \tilde{y}_k)$ converges to $(\xb,\yb)$ as $k\to +\infty$. 

By the definition of sup and inf-convolutions, there exist $x_k \in B(\tilde{x}_k,1)$ and $y_k \in B(\tilde{y}_k,1)$ such that
$$ R^\alpha [u](\tilde{x}_k,p) = u (x_k) - p \cdot (x_k - \tilde{x}_k) 
- \frac{|x_k - \tilde{x}_k|^2}{2 \alpha} \quad \mbox{and}\quad 
 R_\alpha [v](\tilde{y}_k,-q) = v (y_k) - q \cdot (y_k - \tilde{y}_k) 
+ \frac{|y_k - \tilde{y}_k|^2}{2 \alpha}. 
$$
Examining carefully the maximum point property for $(\tilde{x}_k, \tilde{y}_k)$ and using the definition of $ R^\alpha [u],  R_\alpha [v]$, we deduce that
$(x_k,y_k)$ is a maximum point of $(x,y) \mapsto u(x)-v(y)-\phi_k(x,y)$ where $\phi_k (x,y)=\tilde\phi_k(x+\tilde{x}_k-x_k,y+\tilde{y}_k-y_k)$.

We next recall that $R^\alpha [u](\cdot,p)$ and $R_\alpha [v](\cdot,-q)$  are twice differentiable at $\tilde{x}_k$ and $\tilde{y}_k$ respectively.
From property \eqref{eq:inclus} in Proposition~\ref{prop:infconv} together with the a first-order optimality condition for the maximum point 
$(\tilde{x}_k,\tilde{y}_k)$, we deduce that $p_k := \nabla R^\alpha [u](\tilde{x}_k,p)$, 
$X_k := D^2 R^\alpha [u](\tilde{x}_k,p)$, $q_k:= - \nabla R_\alpha [v] (\tilde{y}_k,-q)$, 
$Y_k:= D^2 R_\alpha [v] (y_k,-q)$ satisfy 
\begin{eqnarray*}
(p_k, X_k) \in J^+u(x_k), \quad (-q_k,Y_k) \in J^-v(y_k) \, , \\
 p_k = p - \frac{(\tilde{x}_k-x_k)}{\alpha} = D_x \tilde \phi_k (\tilde{x}_k, \tilde{y}_k)\; , \\
 q_k = q - \frac{(\tilde{y}_k-y_k)}{\alpha}= D_y \tilde \phi_k (\tilde{x}_k, \tilde{y}_k)\; .
\end{eqnarray*}
On the other hand, since $R^\alpha [\phi ] \geq \phi$ and  \eqref{eq:egal-phi} holds, we know that $D R^\alpha [\phi ]((\bar{x},\bar{y}), (p,q))  = D\phi(\bar{x},\bar{y})= (p,q)$.
Recalling that $\tilde \phi_k$ converges in $C^2$ to $R^\alpha [\phi ]((x,y), (p,q))$ and 
 $(\tilde{x}_k, \tilde{y}_k)\to (\xb,\yb)$ as $k\to +\infty$, we deduce that $D_x \tilde \phi_k (\tilde{x}_k, \tilde{y}_k)\to p$ and $D_y \tilde \phi_k (\tilde{x}_k, \tilde{y}_k) \to q$
 as $k\to +\infty$. From the above properties, this yields $\tilde{x}_k-x_k, \tilde{y}_k-y_k \to 0$ and therefore $(x_k,y_k)\to (\xb,\yb)$ as $k\to +\infty$. 

And all the claims of Lemma~\ref{lem:global-jil} are either already proved or of classical properties.
\end{proof}

\subsection{How to apply the lemma?}
\label{subsec:appl-lem}

Now we address the question: how to apply the nonlocal Jensen-Ishii's Lemma? 
The (partial) answer is given by the

\begin{cor}\label{resatten} Let $u$ be an usc viscosity subsolution of \eqref{model},
 let $v$ be a lsc viscosity supersolution of \eqref{model} and let 
$\phi \in  C^2 (\R^d)$. If $(\xb,  \yb) \in \R^{2d}$ is a 
global maximum point of $u(x) -v(y) - \phi (x,y)$, then, for any $\delta >0$, 
there exists $\bar \alpha$ such that, for $0<\alpha < \bar \alpha$, we have
$$ F(\xb,u (\xb) , p , X , \mathcal{I}^{1,\delta}[\xb,\phi_\alpha(\cdot,\yb)]+
 \mathcal{I}^{2,\delta}[\xb,p,u]) \leq 0\; , $$
$$ F(\yb ,v (\yb) , q , Y , \mathcal{I}^{1,\delta}[\yb ,-\phi_\alpha(\xb,\cdot)] +
 \mathcal{I}^{2,\delta}[\yb,q,v]) \geq 0\; ,$$
where $p=\nabla_x\phi (\xb, \yb)=\nabla_x \phi_\alpha (\xb,\yb)$, 
$q = -\nabla_y \phi (\xb, \yb)=\nabla_y \phi_\alpha (\xb,\yb)$ and with also
\begin{equation}\label{ineg-mat-cor}
- \frac{1}{\alpha} I \le \left[\begin{array}{ll} X & 0 \\ 0 & -Y \end{array} \right]
\le D^2\phi_\alpha (\xb, \yb) = D^2\phi (\xb, \yb) + o_\alpha (1)\; .
\end{equation}
\end{cor}

This situation is exactly the one which is needed in the uniqueness proofs, including
 the ones which consist in proving first that $u-v$ is a subsolution of an auxiliary PIDE
(see \cite{bbp97} for details).
 
 It is worth pointing out that, in the nonlocal term, one has a priori to use the function 
$\phi_\alpha$ instead of $\phi$: this is consistent with the fact that the second derivatives 
are also estimated by $D^2 \phi_\alpha (\xb, \yb)$. However, because of our assumption on the
 $\mathcal{I}^{1,\delta}$-term, the difference between 
$\mathcal{I}^{1,\delta}[\xb,\phi_\alpha]$ and $\mathcal{I}^{1,\delta}[\xb,\phi]$ is 
a small error term and we will see in the next section that 
under mild additional assumptions on $F$, one can use 
$\mathcal{I}^{1,\delta}[\xb,\phi] + o_\alpha (1)$ as well.

\begin{proof} Without loss of generality, changing $\phi$ in $\phi -M$ for 
some well-choosen constant $M$, we can assume 
that $(\xb,  \yb) \in \R^{2d}$ is a zero maximum point of $u(x) -v(y) - \phi (x,y)$.

Applying the Lemma together with Definition~\ref{defi:limiting}, 
we get
$$ 
F(x_k,u (x_k) , p_k , X_k , \mathcal{I}^{1,\delta}[x_k,\phi_k(\cdot,y_k)] 
+ \mathcal{I}^{2,\delta}[x_k,\nabla_x \phi_k(x_k,y_k),u])
 \leq 0\; ,
$$
$$ 
F(y_k,v (y_k) , q_k , Y_k , \mathcal{I}^{1,\delta}[y_k,-\phi_k(x_k,\cdot)] 
+ \mathcal{I}^{2,\delta}[y_k,-\nabla_y \phi_k(x_k,y_k),v]) 
\geq 0\; ,
$$
Choosing $\alpha$ small enough in order that $\phi_k \to \phi_\alpha 
: = R^\alpha [\phi]((x,y),(p,q))$ in $C^2(B(x,R_\delta))$, we can pass to the 
limit and obtain the result.
\end{proof}

\medskip

We conclude this section by an easy extension of the nonlocal Jensen-Ishii Lemma. 
It concerns the case of time-dependent equations
\begin{equation}\label{model-par}
 u_t + F(x,u,Du,D^2 u,  \mathcal{I}[x,u]) = 0 \quad\hbox{in  }\R^d \times (0,T)\; ,
\end{equation}
where $T>0$.

We formulate without proof the analogue of Corollary~\ref{resatten} where $J^+ , J^-$ 
denotes the ``parabolic'' super and subjets which take only into account 
the second-order derivatives in space (and not in time).

\begin{cor}\label{resatten-par} Let $u$ be an usc viscosity subsolution 
of \eqref{model-par}, let $v$ be a lsc viscosity supersolution of \eqref{model-par} 
and let $\phi \in  C^2 (\R^d \times (0,T))$. If 
$(\xb,  \yb, \tb) \in \R^{2d}\times (0,T)$ is a global maximum point of 
$u(x,t) -v(y,t) - \phi (x,y,t)$, then, for any $\delta >0$, there exists 
$\bar \alpha$ such that, for $0<\alpha < \bar \alpha$, there exists 
$(a,p,X) \in J^+ u(\xb, \tb)$, $(b,q,Y) \in J^- v(\yb, \tb)$ such that we have
$$ a + F(\xb,u (\xb) , p , X , \mathcal{I}^{1,\delta}[\xb,\phi_\alpha (\cdot,\yb)]+
 \mathcal{I}^{2,\delta}[\xb,p,u]) \leq 0\; , $$
$$ b+ F(\yb ,v (\yb) , q , Y , \mathcal{I}^{1,\delta}[\yb ,-\phi_\alpha(\xb,\cdot)] +
 \mathcal{I}^{2,\delta}[\yb,q,v]) \geq 0\; ,$$
with, in addition, $p=\nabla_x\phi (\xb, \yb)$, $q = -\nabla_y \phi (\xb, \yb)$ and
$$ a-b = \phi_t (\xb, \yb)\; ,$$
$$- \frac{1}{\alpha} I \le \left[\begin{array}{ll} X & 0 \\ 0 & -Y \end{array} \right]
\le D^2\phi_\alpha (\xb, \yb) = D^2\phi (\xb, \yb) + o_\alpha (1)\; .
$$
\end{cor}

The proof of Corollary~\ref{resatten-par} follows readily the classical ideas 
to prove the local Jensen-Ishii Lemma and the above arguments to treat the nonlocal part.

\section{Application to comparison principles}
\label{sec:comparison}

In this section, we consider the equation
\begin{equation}\label{eq:second}
F (x, u , \nabla u, D^2 u , \mathcal{I}_{LI}[u](x) )= 0 \qquad \mbox{~in }\R^d
\end{equation}
where $\mathcal{I}_{LI}[\phi ]$ is given by \eqref{op:levy-ito}. We state and prove a 
comparison principle in the class of bounded (sub  and super) solutions.
We treat a case where there is no strong interaction between the nonlocal term
and the $x$-dependence of $F$; we discuss this assumption in Subsection~\ref{disc-assump}.

\subsection{Statement of the comparison principle}

We first need to strengthen the condition~\eqref{cond:j} on $\mu$ and $j$. 

\begin{itemize}
\item
(A1) The measure $\mu(dz)$ and the function $j(x,z)$ satisfy
\begin{equation}\label{cond:muj1}
\int_{\R^d \setminus B} \mu(dz)<+\infty  \quad \text{ and } \quad 
\sup_{x \in \R^d} \int_B |j(x,z)|^2 \mu(dz) < +\infty \ ,
\end{equation}
and there exists a constant $\bar{c}>0$ such that
\begin{equation}\label{cond:muj2}
 \int_{\R^d
}|j(x,z) - j(y,z)|^2 \mu(dz) \le {\bar c}|x-y|^2\ \mbox{~and }\  \int_{\R^d \setminus B} |j(x,z) -
j(y,z) | \mu (dz)\le {\bar c}|x-y| \ .
\end{equation}
\end{itemize}

For the nonlinearity $F$, we first introduce the classical assumption on the 
dependence of $F$ in $u$.
\begin{itemize}
\item (A2) There exists $\gamma >0$ such that for any $x\in \R^d
$, $u,v\in \R$, $p\in \R^d
$, $X\in \S$ and $l\in \R$
$$ 
 F(x,u,p,X,l)- F(x,v,p, X,l)\ge \gamma (u-v) \qquad \mbox{ when
   } \qquad u \ge v.
$$

\hspace{-\leftmargin}Next we have to impose assumptions on the dependence of $F$ in $x$ 
and we can do it in two ways.

\item (A3-1) For any $R>0$, there exist moduli of continuity $\omega, \omega_R$ such that, 
for any $|x|,|y|\le R$, $|v|\le R$, $l\in\R$ and for any $X,Y \in \S$ satisfying
\begin{equation}\label{ineqmat}
\left[\begin{array}{cc}X&0\\0&-Y\end{array}\right]
 \le \frac1\eps \left[\begin{array}{cc}I&-I\\-I&I\end{array}\right]
+ r(\beta) \left[\begin{array}{cc}I&0\\0&I\end{array}\right]
\end{equation}
for some $\eps >0$ and $r(\beta) \to 0$ as $\beta \to 0$, then, if $s(\beta) \to 0$ 
as $\beta \to 0$, we have
\begin{equation}\label{cond:a31}
F(y,v,\eps^{-1} (x-y), Y,l) - F(x,v,\eps^{-1} (x-y) + s(\beta),X,l)
\le \omega (\beta)+\omega_R (|x-y|+\eps^{-1}|x-y|^2).
\end{equation}
\hspace{-\leftmargin}  or 
\item (A3-2) For any $R>0$, $F$ is uniformly continuous on 
$\R^n \times [_R,R] \times B_R \times D_R \times \R$ where 
$D_R := \{ X \in \S;\ |X| \leq R\}$ and there exist a modulus 
of continuity $\omega_R$ such that, for any $x,y \in \R^d$, 
$|v|\le R$,$ l\in\R$ and for any $X,Y \in \S$ satisfying (\ref{ineqmat}) and
$\eps >0$, we have
\begin{equation}\label{cond:a32}
F(y,v,\eps^{-1} (x-y), Y,l) - F(x,v,\eps^{-1} (x-y),X,l) \le 
\omega_R (\eps^{-1} |x-y|^2 + |x-y|)\; .
\end{equation}
\end{itemize}

\noindent We provide more comments on these assumptions in the next section 
but clearly (A3-1) allows more general dependence in $x$ while (A3-2) allows more general 
dependence in $p$.

\begin{itemize}
\item (A4) $F(x,u,p,X,l)$ is Lipschitz continuous in $l$, uniformly with respect 
to all the other variables. 
\end{itemize}
%-------------- 
\begin{thm}\label{comp}
Assume that (A1), (A2), (A3-1) or (A3-2) and (A4) hold. If $u$ is a bounded usc 
subsolution of \eqref{eq:second} 
and $v$ is a lsc bounded supersolution $v$ of \eqref{eq:second}, then $u \le v$ on $\R^d$. 
\end{thm}
%------
\begin{proof} We consider $M = \sup_{\R^d}(u-v)$ and argue by contradiction by assuming
that $M >0$. 
\medskip

We next approximate $M$ by dedoubling the variables 
$$ 
M_{\eps,\beta}= \sup_{x,y\in \R^d}\left\{u(x)-v(y)-\frac{|x-y|^2}{2\eps}
- \psi_\beta (x)\right\}
$$
where $\eps, \beta$ are small parameters devoted to tend to $0$ and the functions 
$\psi_\beta$ are built in the following way: let $\psi: \R^d \to \R+$ be a smooth 
function such that $\psi, \nabla\psi, D^2 \psi$ are bounded in $\R ^d$ and such that 
$\psi(x) = 0$ if $ |x| \leq 1$ and $\psi(x)>\mathcal{R}:=(||u||_\infty + ||v||_\infty)$ 
for $|x|\ge 2$.
We then set $ \psi_\beta (x) =  \psi (\beta x)$. The three main properties of the 
$\psi_\beta$ we are going to use are the following
\begin{itemize}
\item $ \psi_\beta (x) > \mathcal{R}$ when $|x| \ge 2/\beta$, 
which ensures that the supremum defining $M_{\eps,\beta}$ is achieved and therefore is a 
maximum,
\item $ \nabla \psi_\beta (x), \ D^2 \psi_\beta (x)\to 0$ as $\beta\to 0$ uniformly on 
$\R^d$, which allows to control the differential terms of the $\psi_\beta $,
\item $ \mathcal{I}_{LI}[\psi_\beta](x) \to 0$ as $\beta\to 0$ uniformly on $\R^d$, 
which allows to control the integral terms of the $\psi_\beta $.
\end{itemize}
We classically obtain for $\eps$ and $\beta$ small enough,
\begin{eqnarray*}
0 < \frac{M}2 \le M_{\eps,\beta}\le u(\bar{x})-v(\bar{y})\quad \mbox{and}\quad
\frac{|\bar{x}-\bar{y}|}\eps \le \frac{C}{\sqrt{\eps}}\quad \mbox{and}\quad
\psi_\beta (\bar{x}) \le \mathcal{R}.
\end{eqnarray*}
In particular, $|\bar{x}|\le 2/\beta$. 
Now we consider any maximum points $(\bar{x},\bar{y})$ of the function 
$\displaystyle u(x)-v(y)-\frac{|x-y|^2}{2\eps} - \psi_\beta (x)$. By definition 
of $(\bar{x},\bar{y})$, we have
\begin{eqnarray*}
u (\bar{x}+d')- v (\bar{y} +d)- \frac{|\bar{x}+d'-\bar{y}-d|^2}{2 \eps} 
- \psi_\beta (\bar{x}+d')
\le u (\bar{x}) - v(\bar{y} ) - \frac{|\bar{x}-\bar{y}|^2}{2\eps} -\psi_\beta
 (\bar{x}).
\end{eqnarray*}
By setting $q= \frac{\bar{x}-\bar{y}}\eps$ and $p = q + \nabla \psi_\beta (\bar{x})$, 
we deduce from the previous inequality
\begin{equation}\label{ineq1}
\begin{array}{l}
u (\bar{x}+ j (\bar{x},z))-u (\bar{x}) - p \cdot j (\bar{x},z) \le 
v(\bar{y} + j (\bar{y},z))-v (\bar{y}) - q \cdot j(\bar{y},z)  \\
+ \frac{|j(\bar{x},z)-j(\bar{y},z)|^2}{2\eps}
+ \bigg( \psi_\beta (\bar{x}+j(\bar{x},z))-\psi_\beta(\bar{x})
-\nabla\psi_\beta (\bar{x})\cdot j(\bar{x},z) \bigg) 
\end{array}
\end{equation}
and 
\begin{equation}\label{ineq2}
\begin{array}{l}
u (\bar{x}+ j (\bar{x},z))-u (\bar{x}) \le 
v(\bar{y} + j (\bar{y},z))-v (\bar{y}) + q \cdot (j(\bar{x},z) -
j(\bar{y},z) ) \\
+ \bigg( \psi_\beta (\bar{x}+j(\bar{x},z))-\psi_\beta(\bar{x}) \bigg)
+ \frac{|j(\bar{x},z)-j(\bar{y},z)|^2}{2\eps}.
\end{array}
\end{equation}

For $\delta$ small enough, in order, at least, that $B(0,\delta) \subset B$, 
we define $\mathcal{I}^{1,\delta}$ as being the same integral as $\mathcal{I}_{LI}$ 
but integrating only on $B(0,\delta) $ and, in the same way, $\mathcal{I}^{2,\delta}$ 
stands for the same integral as $\mathcal{I}_{LI}$ but integrating only on 
$\R^d \setminus B(0,\delta)$. We also denote by $\phi_x$ the function 
$x \mapsto \phi(x, \bar{y})$ and by $\phi_y$ the function $y \mapsto \phi(\bar{x},y)$.

Define $\phi (x,y) = \frac{|x-y|^2}{2\eps}+\psi_\beta(x)$. Then an explicit computation of each terms gives
$$ \mathcal{I}^{1,\delta}[\bar{x},\phi_x] = \frac1\eps\int_{|z|\le\delta}|j(\bar{x},z)|^2\mu(dz) 
+ \mathcal{I}^{1,\delta}[\bar{x},\psi_\beta]\; ,$$
$$\mathcal{I}^{1,\delta}[\bar{y},-\phi_y] = - \frac1\eps\int_{|z|\le\delta}|j(\bar{y},z)|^2\mu(dz) \; .$$
Therefore
\begin{eqnarray*}
\mathcal{I}^{1,\delta}[\bar{x},\phi_x] 
\le \mathcal{I}^{1,\delta}[\bar{y},-\phi_y] + 
\frac1\eps\int_{|z|\le\delta}|j(\bar{x},z)|^2\mu(dz) 
+ \frac1\eps\int_{|z|\le\delta}|j(\bar{y},z)|^2\mu(dz) 
+ \mathcal{I}^{1,\delta}[\bar{x},\psi_\beta]
\\
\le \mathcal{I}^{1,\delta}[\bar{y},-\phi_y] +  \frac1\eps o_\delta(1)+o_\beta(1)
\end{eqnarray*}
(we used Assumption (A1)).

Next we consider the $\mathcal{I}^{2,\delta}$ terms which, in fact, consist in two terms,
 whether we integrate on $B \setminus B(0, \delta)$ or on $\R^d \setminus B$. The corresponding 
terms are denoted respectively by
$\mathcal{I}^{2,\delta}_1$ and $\mathcal{I}^{2,\delta}_2$.

For the $\mathcal{I}^{2,\delta}_1$ term, we integrate inequality 
\eqref{ineq1} on $B \cap \setminus B(0, \delta)$, which yields
$$ 
\mathcal{I}^{2,\delta}_1 [\bar{x},p,u]\le 
\mathcal{I}^{2,\delta}_1 [\bar{y},q,v]
+ \int_B \frac{|j(\bar{x},z)-j(\bar{y},z)|^2}{2\eps} \mu (dz)
+\mathcal{I}^{2,\delta}_1 [\bar{x},\nabla\psi_\beta(\bar{x}),\psi_\beta]\; ,
$$
where for the second term of the right-hand side, we have estimated 
the integral on $B \cap ( \R^d \setminus B(0, \delta))$ by the integral on the whole ball $B$.

Next, for the $\mathcal{I}^{2,\delta}_2$ term, we integrate inequality 
\eqref{ineq2} on $\R^d \setminus B$, which yields
\begin{eqnarray*}
 \mathcal{I}^{2,\delta}_2 [\bar{x},u]\le 
\mathcal{I}^{2,\delta}_2 [\bar{y},v]
+ \int_{\R^d \setminus B} q \cdot (j(\bar{x},z) -j(\bar{y},z) ) \mu (dz)\\
+ \frac1{2\eps} \int_{\R^d \setminus B} |j(\bar{x},z)-j(\bar{y},z)|^2 \mu (dz)
+\mathcal{I}^{2,\delta}_2 [\bar{x},\nabla\psi_\beta(\bar{x}),\psi_\beta]
\; .
\end{eqnarray*}

By using \eqref{cond:muj1}-\eqref{cond:muj2}, summing up these inequalities, we thus obtain
$$ 
\mathcal{I}^{2,\delta} [\bar{x},p,u]\le 
\mathcal{I}^{2,\delta} [\bar{y},q,v]
+ O\left(\frac{|\bar{x}-\bar{y}|^2}{\eps}\right) + o_\beta (1)\; 
$$
and finally, we get the following estimate between the integral terms
\begin{equation}\label{estim:int}
l:= \mathcal{I}^{1,\delta}[\bar{x},\phi_x]
+\mathcal{I}^{2,\delta}[\bar{x},p,u] \leq \mathcal{I}^{1,\delta}[\bar{y},-\phi_y]
+\mathcal{I}^{2,\delta}[\bar{y},q,v] +O\left(\frac{|\bar{x}-\bar{y}|^2}{\eps}\right) 
+ o_\beta (1)  + \frac1\eps o_\delta (1)\; .
\end{equation}
\medskip

We are next going to apply Corollary~\ref{resatten} with $\phi$. 
If $(p,-q)$ denotes $\nabla \phi (\bar{x},\bar{y})$ and  $A$ denotes 
$D^2 \phi (\bar{x},\bar{y})$, for any $\alpha >0$, there exists 
two matrices $X,Y \in \S$ such that \eqref{ineg-mat-cor} holds
true and such that
\begin{eqnarray*}
F (\bar{x}, u (\bar{x}),p,X,\mathcal{I}^{1,\delta}[\bar{x},\phi_\alpha(\cdot, \bar{y})]
+\mathcal{I}^{2,\delta}[\bar{x},p,u])\le 0 \ , \\
F(\bar{y}, v (\bar{y}),q,Y,\mathcal{I}^{1,\delta}[\bar{y},-\phi_\alpha (\bar{x},\cdot)]
+\mathcal{I}^{2,\delta}[\bar{y},q,v]) \ge 0 
\end{eqnarray*}
where $\phi_\alpha$ is defined in Lemma~\ref{lem:global-jil}.
We next use Proposition~\ref{prop:infconv} in order to get
\begin{eqnarray*}
F (\bar{x}, u (\bar{x}),p,X,l)\le o_\alpha (1) \ , \\
F(\bar{y}, v (\bar{y}),q,Y,\mathcal{I}^{1,\delta}[\bar{y},-\phi_y]
+\mathcal{I}^{2,\delta}[\bar{y},q,v]) \ge o_\alpha (1) \, .
\end{eqnarray*}
We next combine the two previous viscosity inequalities and we use  (A2)
and (A4) together with Estimate~\eqref{estim:int} and the (nonlocal) ellipticity 
in order to get 
\begin{equation}\label{estim:visc}
\gamma\frac{M}2 \le F(\bar{y},v(\bar{y}),q,Y,l)
-F(\bar{x},v(\bar{y}),q+\nabla \psi_\beta(\bar{x}),X,l) 
+O\left(\frac{|\bar{x}-\bar{y}|^2}{\eps}\right) 
+ o_\beta (1) + o_\alpha (1) + \frac1\eps o_\delta (1)\ .
\end{equation}
Inequality~\eqref{ineg-mat-cor} implies in particular that
\begin{eqnarray*}
- \frac1\alpha \left[\begin{array}{cc}I&0\\0&I\end{array}\right]\le
\left[\begin{array}{cc}X&0\\0&-Y\end{array}\right]
 \le \frac1{\eps}\left[\begin{array}{cc}I&-I\\-I&I\end{array}\right]+
\left[\begin{array}{cc} D^2 \psi_\beta (\bar{x})&0\\0&0\end{array}\right] 
+ o_\alpha (1)\left[\begin{array}{cc}I&0\\0&I\end{array}\right]
\\
  \le \frac1{\eps}\left[\begin{array}{cc}I&-I\\-I&I\end{array}\right]
+ (o_\alpha(1)+o_\beta (1))\left[\begin{array}{cc}I&0\\0&I\end{array}\right]
\end{eqnarray*}
where $o_\beta (1)$ and $o_\alpha (1)$ are uniform in $\eps$;
 we used the properties satisfied by $\psi_\beta$ listed above. 
%In the following, we will choose $\alpha=1$ so that 
%$|X|+|Y|\le {\bar c}(\eps^{-1}+o_\beta (1)+o_\eps (1))$.
\medskip

We first assume that (A3-1) holds. In this case, we choose $R=R_\beta=2/\beta$ and write 
\begin{multline*}
\gamma\frac{M}2 \le  F(\bar{x},v(\bar{y}),q,Y,l)
-F(\bar{y},v(\bar{y}),q+o_\beta(1),X,l) 
+O\left(\frac{|\bar{x}-\bar{y}|^2}{\eps}\right) + o_\beta (1)  + o_\alpha (1) 
+ \frac1\eps o_\delta (1)\\
\le  \omega \left( \frac{|\bar{x}-\bar{y}|^2}\eps + |\bar{x}-\bar{y}|
+ o_\beta (1) \right)+\omega_{R_\beta} \left( \frac{|\bar{x}-\bar{y}|^2}\eps + |\bar{x}-\bar{y}|
 \right)+O\left(\frac{|\bar{x}-\bar{y}|^2}{\eps}\right) + o_\beta (1) + o_\alpha (1) 
+ \frac1\eps o_\delta (1) \ .
\end{multline*}
Using \eqref{cond:a31} and letting successively $\delta$, $\alpha$, $\eps$ 
and $\beta$ go to $0$, we get the desired result.
\medskip

Next we assume that (A3-2) holds. We derive from \eqref{estim:visc}, 
\begin{eqnarray*}
\gamma\frac{M}2 \le F(\bar{x},v(\bar{y}),q,Y,l)
-F(\bar{y},v(\bar{y}),q,X,l) 
+O\left(\frac{|\bar{x}-\bar{y}|^2}{\eps}\right) 
+ \omega_{R_\eps} (|\nabla\psi_\beta|)  + o_\beta (1) 
+ o_\alpha (1) + \frac1\eps o_\delta (1)\\
\le \omega \left(\frac{|\bar{x}-\bar{y}|^2}{\eps}+|\bar{x}-\bar{y}|+o_\beta(1)\right)
+O\left(\frac{|\bar{x}-\bar{y}|^2}{\eps}\right) 
+  \omega_{R_\eps} (|\nabla\psi_\beta|) + o_\beta (1) + o_\alpha (1) + \frac1\eps o_\delta (1)
\end{eqnarray*}
with $R_\eps = K/{\eps}$ for some $K>0$. In this case, we use \eqref{cond:a32} and we let successively $\delta$, $\alpha$, 
$\beta$ and $\eps$ tend to $0$. The proof is now complete. 
\end{proof}

\begin{rem} It is worth pointing out that, in the two cases we consider in the proof (letting first $\eps$ tend to $0$ and then $\beta$ tend to $0$, or the contrary), the behavior of $M_{\eps, \beta}$ are different. Indeed 
$$  \lim_{\beta \to 0} \lim_{\eps \to 0} M_{\eps, \beta} =M\quad \hbox{while }\quad
 \lim_{\eps \to 0} \lim_{\beta \to 0} M_{\eps, \beta} = \lim_{s \downarrow 0} \sup_{|x-y|\le s} (u(x)-v(y)) \geq M\; .$$ But, in both cases, we have the key property $\displaystyle \frac{|\bar{x}-\bar{y}|^2}{\eps} \to 0$.
\end{rem}

\subsection{Discussion of the assumptions}\label{disc-assump}

We want to discuss here the assumptions (A3-1) and (A3-2) and, in particular, the connections with the nonlocal term.

First, as long as local equations are concerned, we recall that the equation
$$ - {\rm Tr} (A(x) D^2u) + b(x)|Du|^p + c(x) u = f(x) \quad \hbox{in  }\R^n\; ,$$
satisfies (A3-1) if $A= \sigma^T \sigma$ where $\sigma$ is a bounded, matrix-valued locally Lipschitz continuous function, $0 < p \leq 1$ and $b$ is also a locally Lipschitz continuous function and $c, f$ are continuous functions. For (A3-2), in most of the cases, we have still to assume $0 < p \leq 1$ but $\sigma$ and $b$ have to be (globally) bounded Lipschitz continuous functions and $c,f$ need to be uniformly continuous. But if $b$ is a constant function then $p$ can be any nonnegative number.

It is worth pointing out that the assumptions on the nonlocal term are rather restrictive: for example, one can add (in fact subtract) the L\'evy operator $\mathcal{I}_L [u](x)$ given by \eqref{op:levy} from the above equation and the resulting equation still satisfies (A3-1) or (A3-2). But this is not the case anymore if the measure $\mu(dz)$ has a singularity at $z=0$ and if we try to subtract a term like $d(x)\mathcal{I}_L [u](x)$, whatever we may assume on the function $d$. We can treat such term only if $\mu(dz)$ is a bounded measure. Therefore the $x$-dependence in the nonlocal term is rather restrictive, except perhaps if it written in the L\'evy-Ito form \eqref{op:levy-ito} where we have a well-adapted dependence in $x$. This is the reason why we formulate Theorem~\ref{comp} with such an operator.

Curiously the type of singularity of $\mu(dz)$ does not seem to change anything: one could think that the cases where $|z|$ is integrable at $0$ and the cases where only $|z|^2$ is integrable are different and the first one easier to treat. But we were unable to see any difference in the proof where we can just just use inequalities \eqref{ineq1} and \eqref{ineq2}.

%-------------------------------------------------------------------------

\end{document}